# HYBRID SHRINKAGE ESTIMATORS USING PENALTY BASES FOR THE ORDINAL ONE-WAY LAYOUT

By Rudolf Beran[1]

*University of California, Davis*

This paper constructs improved estimators of the means in the Gaussian saturated one-way layout with an ordinal factor. The least squares estimator for the mean vector in this saturated model is usually inadmissible. The hybrid shrinkage estimators of this paper exploit the possibility of slow variation in the dependence of the means on the ordered factor levels but do not assume it and respond well to faster variation if present. To motivate the development, candidate penalized least squares (PLS) estimators for the mean vector of a one-way layout are represented as shrinkage estimators relative to the penalty basis for the regression space. This canonical representation suggests further classes of candidate estimators for the unknown means: monotone shrinkage (MS) estimators or soft-thresholding (ST) estimators or, most generally, hybrid shrinkage (HS) estimators that combine the preceding two strategies. Adaptation selects the estimator within a candidate class that minimizes estimated risk. Under the Gaussian saturated one-way layout model, such adaptive estimators minimize risk asymptotically over the class of candidate estimators as the number of factor levels tends to infinity. Thereby, adaptive HS estimators asymptotically dominate adaptive MS and adaptive ST estimators as well as the least squares estimator. Local annihilators of polynomials, among them difference operators, generate penalty bases suitable for a range of numerical examples. In case studies, adaptive HS estimators recover high frequency details in the mean vector more reliably than PLS or MS estimators and low frequency details more reliably than ST estimators.

**1. Introduction.** Consider the one-way layout of ANOVA. A single factor that influences the observed responses has $p$ distinct levels $\{s_i : 1 \leq i \leq p\}$. These factor levels can be either nominal (i.e., pure labels that bear no

Received August 2001; revised May 2003.
[1]Supported in part by NSF Grant DMS-03-00806.

ordering information) or ordinal (i.e., real numbers whose order and spacing carries information). In the case of an ordinal factor, we will suppose that the factor levels have been ordered from smallest to largest. At level $s_i$, we observe measurements $\{y_{ij} : 1 \leq j \leq n_i\}$. The *saturated* Gaussian model for the one-way layout asserts that the observations $\{y_{ij}\}$ satisfy

$$(1.1) \qquad y_{ij} = \mu_i + e_{ij}, \qquad 1 \leq i \leq p,\, 1 \leq j \leq n_i.$$

Here the errors $e = \{e_{ij}\}$ are independent, identically distributed, each having an $N(0, \sigma^2)$ distribution and the means $\{\mu_i\}$ are unknown real numbers subject to *no* restrictions. That the means depend on the respective factor levels can be expressed formally by

$$(1.2) \qquad \mu_i = m(s_i), \qquad 1 \leq i \leq p.$$

In equation (1.2), the function $m$ is real-valued, unknown, and is subject to *no* restrictions.

At first glance, the saturated one-way layout model expressed by equations (1.1) and (1.2) resembles a model for curve estimation. However, there is a fundamental distinction. In curve estimation, the domain of $m$ is a continuum, usually a closed subset of the real line. In the one-way layout, the domain of the function $m$ is a discrete set of factor levels. Even in ordinal one-way layouts, no credible extension of $m$ to a larger domain may exist. Tukey [(1977), Chapter 7] fitted several examples of ordinal one-way layouts that are not curve estimation problems because of intrinsic limitations on the domain of the function $m$.

*Hereafter, unless otherwise stated, we consider only ordinal one-way layouts.* The following examples will serve as test cases for our methods:

EXAMPLE 1. The top subplot in Figure 1 displays monthly Australian red wine sales (in kiloliters) from January 1980 to October 1991. The data was reported by Brockwell and Davis (1996) and was analyzed there with techniques based on ARMA models. ARMA models are only one class of hypothetical probability models that might be entertained as a way of mimicking the wine sales data. Because the data is not actually random, it is prudent to carry out alternative analyses. As Tukey (1980) pointed out, "In practice, methodologies have no assumptions and deliver no certainties." We will analyze the wine-sales data with mean estimators derived for the ordinal one-way layout model. Motivating this approach is the traditional decomposition of an econometric times series into a deterministic term (trend plus seasonal variation), plus a random noise term. The factor levels are the 142 successive months in the period considered and are clearly ordinal. Ipso facto, mean monthly wine-sales are defined only on the discrete time grid of months. Our analysis in Section 2.5 finds a highly intelligible seasonal pattern in the wine sales.

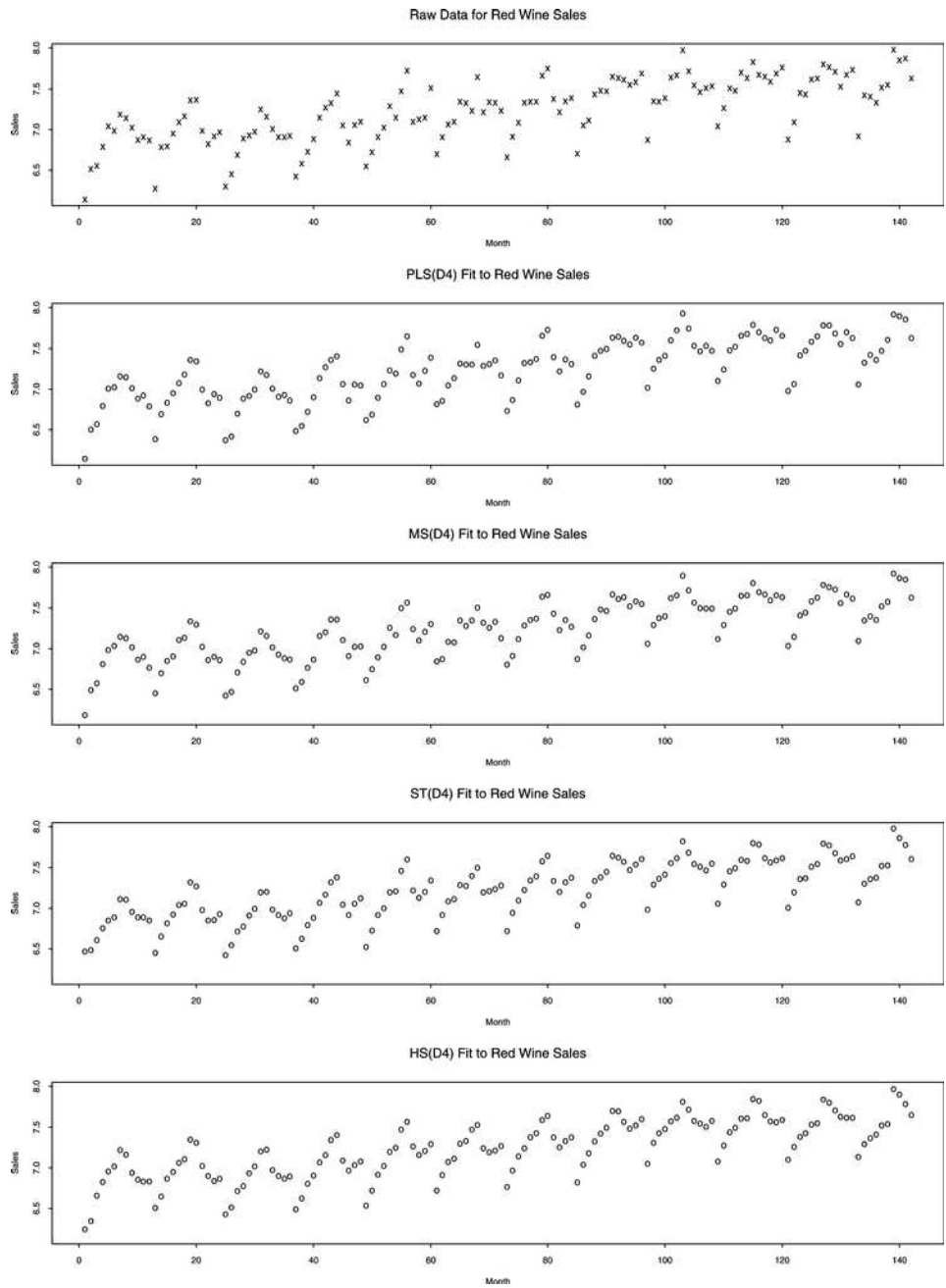

Fig. 1. *Competing* D4-*basis fits to the Australian monthly red wine-sales data.*

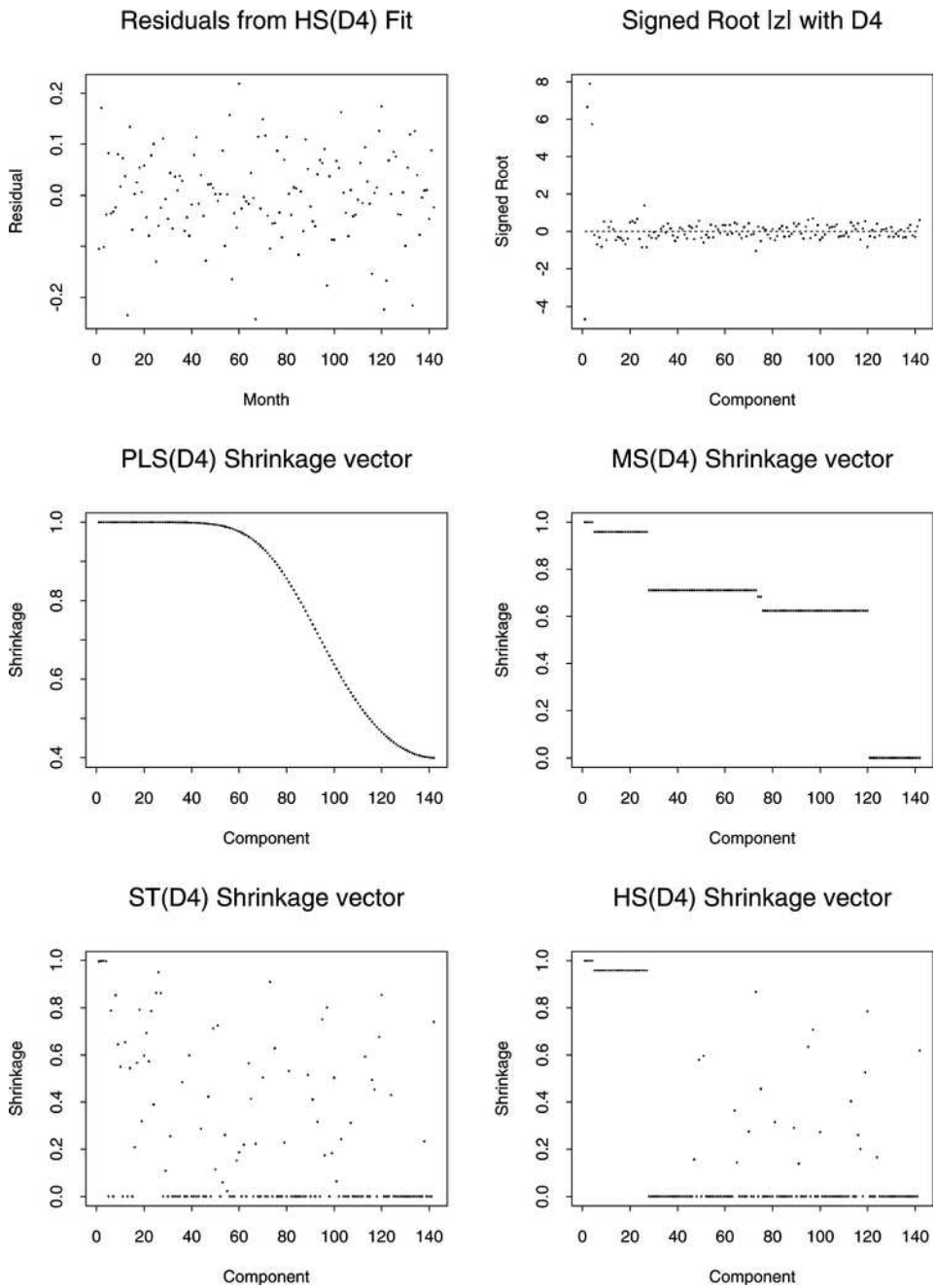

FIG. 2. *Diagnostics for* D4-*basis fits to the Australian monthly red wine-sales data: residuals for the* HS(D4) *fit, the empirical basis economy plot and the shrinkage vectors used by competing* D4-*basis estimators.*

EXAMPLE 2. The artificial ordinal one-way layouts analyzed in Figures 3 and 4 are designed to bracket the situation found in the case study of Example 1. In each of Figures 3 and 4, the data in the top subplot is obtained by adding pseudo-random errors to the means displayed in the second subplot. The means in Figure 3 vary slowly while those in Figure 4 vary rapidly. To the human eye, the pattern of variation in the means is not visible in the data. In Section 2.6, comparing competing estimators of means on these two artificial ordinal one-way layouts adds to our understanding of their performance.

Form the $n \times 1$ observation vector $y = \{\{y_{ij} : 1 \leq j \leq n_i\}, 1 \leq i \leq p\}$, where $n$ is the total number of observations. Let $X$ be the $n \times p$ incidence matrix that links observations to the relevant factor level. The $i$th column of $X$ contains $n_i$ ones, the other elements being zeroes. Let $\mu = (\mu_1, \mu_2, \ldots, \mu_p)'$, where $\mu_i$ satisfies (1.2) with $m$ unrestricted. The saturated model (1.1) is equivalent to the assertion

$$(1.3) \qquad y \sim N(\eta, \sigma^2 I_n) \qquad \text{where } \eta = X\mu.$$

The primary task in this paper is to devise regularized estimators of $\eta$, or, equivalently, of $\mu = (X'X)^{-1} X' \eta$, that (asymptotically in $p$) dominate the least squares estimator $\hat{\eta}_{\text{LS}} = X(X'X)^{-1} X' y$ under the saturated ordinal model. We note that the desirability of analyzing the risk of estimators of $\eta$ under the saturated model is a basic way in which estimation in the one-way layout differs from curve estimation.

Suppose that we assess any estimator $\hat{\eta}$ through its normalized quadratic loss and corresponding risk

$$(1.4) \qquad L(\hat{\eta}, \mu) = p^{-1} |\hat{\eta} - \eta|^2, \qquad R(\hat{\eta}, \eta, \sigma^2) = \mathrm{E} L(\hat{\eta}, \eta),$$

the expectation being calculated under the saturated model. Equivalently, we could discuss estimation of $\mu$ under the loss function $p^{-1}(\hat{\mu} - \mu)' X'X (\hat{\mu} - \mu)$. The risk of $\hat{\eta}_{\text{LS}}$ is evidently $\sigma^2$. It is well known that this value is the smallest risk attainable by unbiased estimators of $\eta$ in the saturated model whether the factor is nominal or ordinal. Nevertheless, for both types of factor, $\hat{\eta}_{\text{LS}}$ is an inadmissible estimator of $\eta$ whenever the number $p$ of factor levels exceeds two [Stein (1956)].

The James–Stein (1961) shrinkage estimator of $\eta$ improves significantly on the quadratic risk of $\hat{\eta}_{\text{LS}}$ and is a good answer when the factor is nominal. For an ordinal factor, estimators for $\eta$ that have still lower risk in the one-way layout are often possible. The better estimators of $\eta$ developed in this paper rely on a regularization strategy that enables the data to influence estimator construction. Our hybrid shrinkage estimators exploit the possibility of slow variation in the dependence of the means on the ordered factor levels, but do *not* assume it, and respond well to faster variation if present.

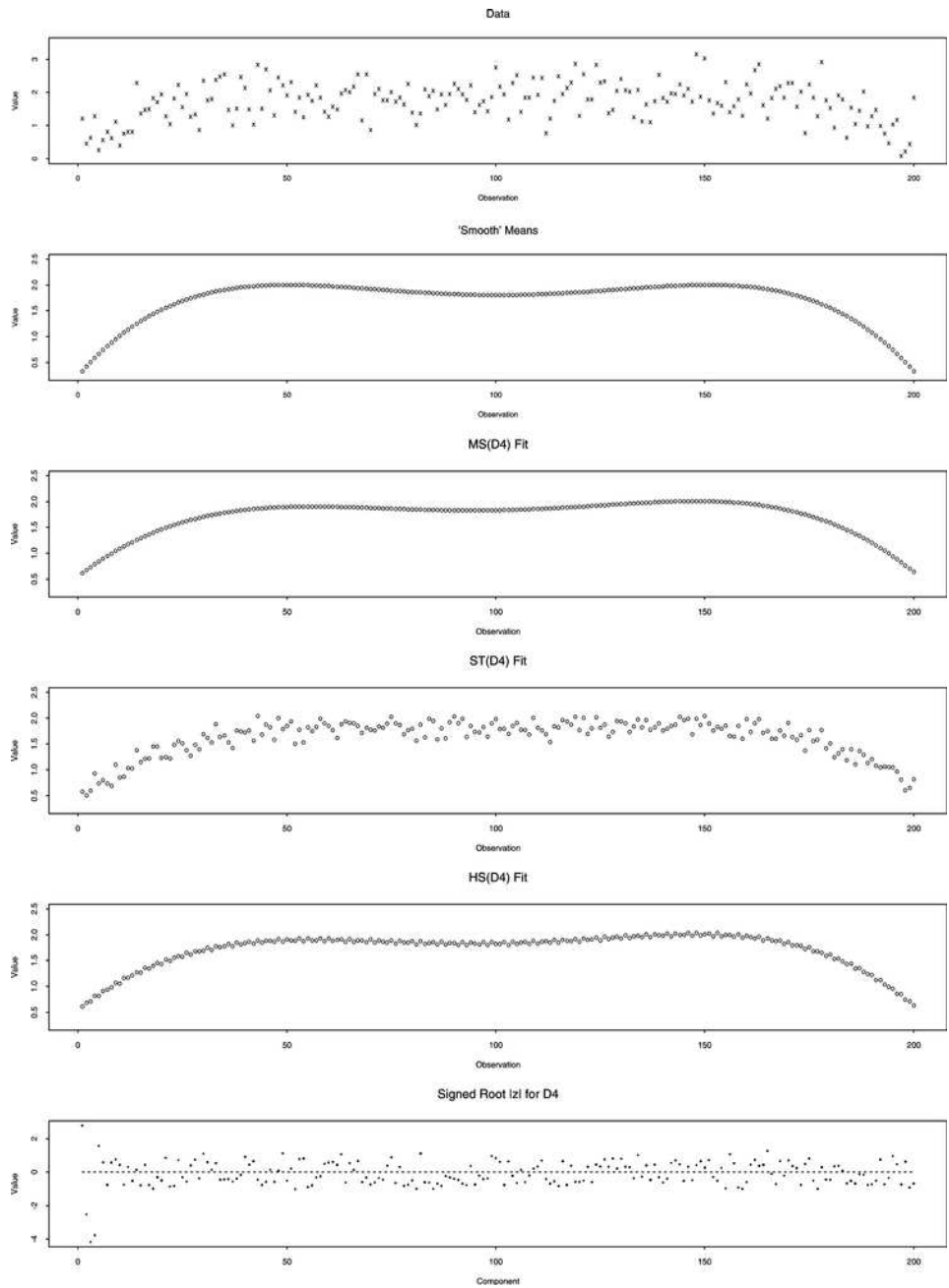

Fig. 3. *Competing D4-basis fits to the Smooth artificial data and the empirical basis economy plot.*

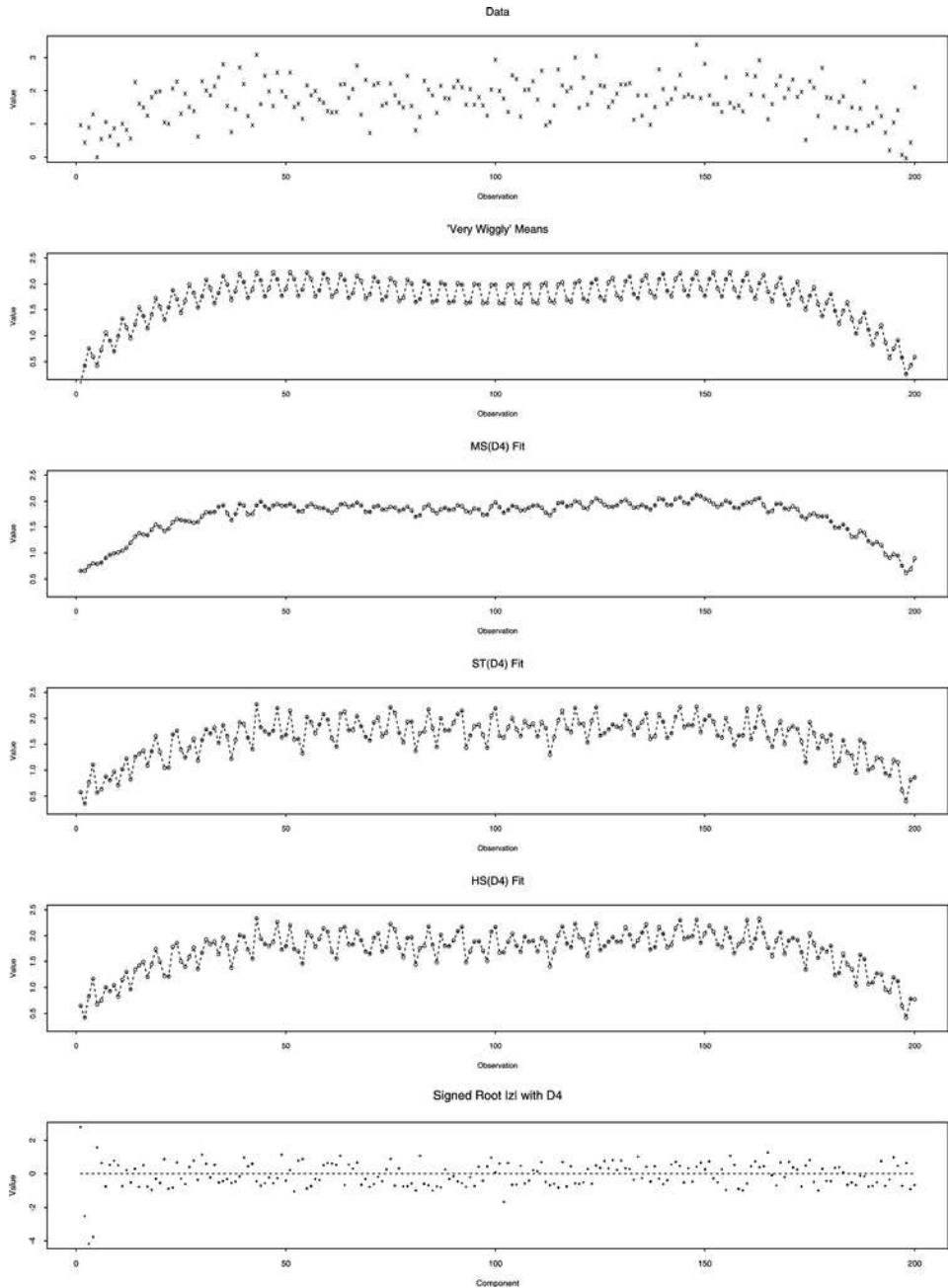

FIG. 4. *Competing* D4-*basis fits to the Very Wiggly artificial data and the empirical basis economy plot. Interpolating lines are added to guide the eye through the sequence of means or estimated means. They have no further significance.*

The broad approach is the following: (a) use prior conjecture about the unknown means in the Gaussian saturated one-way layout to motivate classes of candidate estimators for these means; (b) estimate the risk of each candidate estimator under the saturated model; (c) define an adaptive estimator to be the candidate procedure with smallest estimated risk; (d) experiment with the adaptive estimator on both observed and artificial data; (e) study the asymptotic risk of such adaptive estimators under the saturated model.

The inadmissibility of least squares fits to the means of a Gaussian saturated one-way layout has inspired considerable work on competing estimators. Candidate model selection, ridge regression or penalized least squares (PLS) estimators are all particular symmetric linear estimators. Important studies of symmetric linear estimators include Stein (1981), Li and Hwang (1984), Buja, Hastie and Tibshirani (1989) and Kneip (1994). Tukey (1977) proposed and experimented with certain smoothing algorithms for fitting ordinal one-way layouts. Beran and Dümbgen (1998) used a finite-dimensional version of Pinsker's (1980) asymptotic minimax bound to assess adaptive symmetric linear estimators that perform monotone shrinkage relative to a fixed orthonormal basis.

Adaptive hybrid shrinkage (HS) estimators for the vector $\eta$, the main contribution of this paper, combine monotone shrinkage (MS)—a generalization of PLS—with the soft-thresholding (ST) idea in Donoho and Johnstone (1995). The adaptive HS estimators are devised to dominate asymptotically both adaptive MS and adaptive ST estimators of $\eta$. Theorem 4.1 gives the supporting risk analysis under the saturated model as the number $p$ of factor levels tends to infinity. Interpretation of asymptotic minimax Theorem 3.1 isolates basis economy as a key factor in superior performance of MS estimators and approximate basis economy as a key factor in superior performance of HS estimators. Applied to the penalty bases used in this paper, this interpretation suggests that HS estimators behave like MS estimators when the means of an ordinal one-way layout vary slowly and share the superior ability of ST estimators to track means that vary more rapidly. Related to HS estimators in strategy but not in tactics are the hybrid wavelet fits of Efromovich (1999). These combine a certain linear shrinkage strategy with hard-thresholding of wavelet coefficients.

Sections 2.5 and 2.6 continue the analysis of Examples 1 and 2. Comparisons through estimated risks are supplemented by basis economy plots and shrinkage vector plots that reveal working details of the competing estimators. The diagnostic plots in these examples support the claim made above that basis economy is important for superior performance of MS estimators and that approximate basis economy is important for superior performance of HS estimators. In particular, the numerical experiments confirm the superior ability of adaptive HS estimators constructed on $d$th difference penalty bases to recover both low and high frequency features in the means of an ordinal one-way layout.

Curve estimation can be split conceptually into two problems: (a) estimation of means on the ordinal one-way layout of observed factor levels; and (b) estimation of the mean function between adjacent factor levels through some form of interpolation. The choice of function class in curve estimation strongly affects the implicit interpolation scheme. For nonparametric curve estimation, adaptive curve estimators that achieve the Pinsker asymptotic minimax bound over specified function classes were developed by Efromovich and Pinsker (1984) and by Golubev (1987). On the other hand, data does not come with an attached probability model. A data analyst interested in curve estimation, but not certain of an appropriate function class, might reasonably use the techniques of this paper to estimate the means at the observed factor levels; and might then experiment with curve estimates obtained from these by various interpolation schemes.

This paper distinguishes strictly among data, statistical procedure, probability model and pseudo-random numbers. Modern computing environments for applied and experimental statistics have returned the distinctions to prominence. An adaptive procedure implicitly fits the probability model that motivates it. However, using such a procedure on data differs from believing that a probability model governs the data. Data is not certifiably random. Mathematical study of a statistical procedure under a probability model tests the procedure only on virtual data governed by that model. Such mathematical explorations become pertinent to statistical theory if the probability model can approximate salient relative frequencies in actual data of interest. Our understanding of statistical procedures is ultimately empirical, aided considerably by suitable diagnostic plots, knowledge of the substantive field, and intuitive interpretations of relevant mathematical results [cf. Brillinger and Tukey (1985), Section 17, Beran (2001), Section 3, and Friedman (2001)]. In such respects, statistics does not differ from other sciences that address the world around us.

**2. HS estimators.** This section begins by defining PLS estimators for the mean vector of the saturated ordinal one-way layout and then MS or ST estimators that use the same penalty basis. This background enables the definition of HS estimators that combine the MS and ST shrinkage strategies. Adaptive HS estimators are designed to perform well whether the components of the mean vector vary slowly or more rapidly. Our treatment covers both balanced and unbalanced one-way layouts. Section 4 develops asymptotic theory under the saturated model that supports the adaptation methodology used.

2.1. *Canonical representation of PLS estimators.* As described in the Introduction, the saturated model for the ordinal one-way layout with $p$ factor levels asserts that the observation vector $y$ has an $N(\eta, \sigma^2 I_n)$ distribution, where $\eta = X\mu$. Here $X$ is the incidence matrix that links observations to the relevant factor levels and $n$ is the total sample size. The task is to estimate the mean vector $\eta$. Let $D$ be any matrix with $p$ columns, let $\nu$ be an element of the extended nonnegative reals $[0, \infty]$, and let $|\cdot|$ denote quadratic norm. The candidate PLS estimator of $\eta$ is

$$(2.1) \qquad \hat{\eta}_{\text{PLS}}(D, \nu) = X\hat{\mu}_{\text{PLS}}(D, \nu),$$

where

$$(2.2) \qquad \hat{\mu}_{\text{PLS}}(D, \nu) = \arg\min_{\mu \in R^p}[|y - X\mu|^2 + \nu|D\mu|^2].$$

It is understood that $\hat{\mu}_{\text{PLS}}(D, \infty) = \lim_{\nu \to \infty} \hat{\mu}_{\text{PLS}}(D, \nu)$. The foregoing displays yield the explicit formula

$$(2.3) \qquad \hat{\eta}_{\text{PLS}}(D, \nu) = X(X'X + \nu D'D)^{-1}X'y.$$

Both $D$ and $\nu$ are to be chosen so to control the quadratic risk of the PLS estimator under the saturated model.

In the leading case of a balanced one-way layout, the matrix $X'X$ is a multiple of the identity matrix. Consequently, $\hat{\eta}_{\text{PLS}}$ may be computed equivalently by applying the PLS strategy to the averages $\{y_i : 1 \le i \le p\}$, rather than to the original data. Thus, the case $n = p$ implicitly includes the general balanced one-way layout. Of course, estimating $\sigma^2$ is easier when $n$ exceeds $p$ (see Section 2.2).

A revealing canonical representation of $\hat{\eta}_{\text{PLS}}(D, \nu)$ is obtained through the following algebraic reduction. The replication matrix $R = X'X$ is a $p \times p$ diagonal matrix whose $k$th diagonal element is the number of observations at factor level $s_k$. Let $\mathcal{M}$ denote the the regression space of the one-way layout—the subspace spanned by the columns of the incidence matrix $X$. The columns of the matrix $U_0 = XR^{-1/2}$ provide an orthonormal basis for this regression space. Let $B = R^{-1/2}D'DR^{-1/2}$ have spectral representation $B = \Gamma\Lambda\Gamma'$, where the eigenvector matrix satisfies $\Gamma'\Gamma = \Gamma\Gamma' = I_p$ and the diagonal matrix $\Lambda = \text{diag}\{\lambda_i\}$ gives the ordered eigenvalues with $0 \le \lambda_1 \le \lambda_2 \le \cdots \le \lambda_p$. This eigenvalue ordering, the reverse of the customary, is used here because the eigenvectors associated with the smallest eigenvalues largely determine the value and performance of candidate estimator $\hat{\eta}_{\text{PLS}}(D, \nu)$. Let $U = U_0\Gamma$. It follows from (2.3) that

$$(2.4) \qquad \hat{\eta}_{\text{PLS}}(D, \nu) = U(I_p + \nu\Lambda)^{-1}U'y.$$

The columns of the matrix $U$ define the orthonormal *penalty basis* for the regression space $\mathcal{M}$ of the one-way layout.

Let $z = U'y$ and let $f(\nu)$ denote the column vector $(1/(1+\nu\lambda_1), 1/(1+\nu\lambda_2), \ldots, 1/(1+\nu\lambda_p))'$, with the understanding that $f(\infty) = \lim_{\nu \to \infty} f(\nu)$. The distribution of $z$ is then $N_p(\xi, \sigma^2 I_p)$, where $\xi = U'\eta$. The candidate PLS estimator of $\xi$ implied by expression (2.4) is

$$(2.5) \qquad \hat{\xi}_{\text{PLS}}(D, \nu) = U'\hat{\eta}_{\text{PLS}}(D, \nu) = f(\nu)z,$$

where the multiplication of vectors in the expression to the right is performed componentwise as in the S language. Equivalently,

$$(2.6) \qquad \hat{\eta}_{\text{PLS}}(D, \nu) = U\hat{\xi}_{\text{PLS}}(D, \nu) = U \operatorname{diag}\{f(\nu)\} U'y.$$

REMARK. The successive columns $\{u_j : 1 \leq j \leq p\}$ of the penalty basis matrix $U = U_0 \Gamma$, where $U_0 = XR^{-1/2}$, have a variational characterization:

- Let $\gamma_j$ denote the $j$th column of the eigenvector matrix $\Gamma$.
- Find a unit vector $u_1$ in $\mathcal{M}$ that minimizes the penalty $|D(X'X)^{-1}X'u_1|^2$. The answer is $u_1 = U_0\gamma$, where $\gamma$ is a $p \times 1$ unit vector that minimizes $|DR^{-1/2}\gamma|^2 = \gamma'B\gamma$. Thus, $u_1 = U_0\gamma_1$.
- Find a unit vector $u_2$ in $\mathcal{M}$ that minimizes the penalty $|D(X'X)^{-1}X'u_2|^2$ subject to the constraint that $u_2$ is orthogonal to $u_1$. The answer is $u_2 = U\gamma$, where $\gamma$ is a $p \times 1$ unit vector orthogonal to $\gamma_1$ that minimizes $|DR^{-1/2}\gamma|^2 = \gamma'B\gamma$. Thus, $u_2 = U_0\gamma_2$.
- Continue sequential constrained minimization to obtain the penalty basis matrix

$$(2.7) \qquad U = (U_0\gamma_1, U_0\gamma_2, \ldots, U_0\gamma_p) = U_0\Gamma.$$

2.2. *Adaptive MS estimators.* The canonical representation (2.6) of PLS estimators suggests a larger class of candidate shrinkage estimators that use the same penalty basis $U$. Let

$$(2.8) \qquad \mathcal{F}_{\text{MS}} = \mathcal{F}_{\text{MS}}(p) = \{f \in [0,1]^p : f_1 \geq f_2 \geq \cdots \geq f_p\}$$

and let

$$(2.9) \qquad \hat{\xi}_{\text{MS}}(D, f) = fz, \qquad f \in \mathcal{F}_{\text{MS}}.$$

The candidate MS estimators for $\eta$ associated with penalty matrix $D$ are defined by

$$(2.10) \qquad \hat{\eta}_{\text{MS}}(D, f) = U\hat{\xi}_{\text{MS}}(D, f) = U \operatorname{diag}\{f\} U'y, \qquad f \in \mathcal{F}_{\text{MS}}.$$

It follows from (2.6) that the candidate PLS estimators are a proper subset of the MS family in which the shrinkage vector $f$ is restricted to the form $\{f(\nu) : \nu \in [0, \infty]\}$.

For any vector $x$, let $\mathrm{ave}(x)$ denote the average of its components. Define the function

$$r_{\mathrm{MS}}(f, \xi, \sigma^2) = \mathrm{ave}[f^2 \sigma^2 + (1-f)^2 \xi^2], \qquad f \in [0,1]^p. \tag{2.11}$$

Because $|\hat{\eta}_{\mathrm{MS}}(D, f) - \eta|^2 = |fz - \xi|^2$, it follows that the normalized quadratic risk of the candidate MS estimator is

$$R(\hat{\eta}_{\mathrm{MS}}(D, f), \eta, \sigma^2) = r_{\mathrm{MS}}(f, \xi, \sigma^2), \qquad f \in \mathcal{F}_{\mathrm{MS}}. \tag{2.12}$$

In particular, the risk of the candidate PLS estimator is just $r_{\mathrm{MS}}(f(\nu), \xi, \sigma^2)$.

The risk function $r_{\mathrm{MS}}(f, \xi, \sigma^2)$ depends on the unknown parameters $\xi^2$ and $\sigma^2$. Having obtained a variance estimator $\hat{\sigma}^2$, we may estimate $\xi^2$ by $z^2 - \hat{\sigma}^2$ and, hence, the risk function by

$$\begin{aligned}
\hat{r}_{\mathrm{MS}}(D, f) &= \mathrm{ave}[f^2 \hat{\sigma}^2 + (1-f)^2 (z^2 - \hat{\sigma}^2)] \\
&= \mathrm{ave}[(f - \hat{g})^2 z^2] + \hat{\sigma}^2 \mathrm{ave}(\hat{g}),
\end{aligned} \tag{2.13}$$

where $f \in \mathcal{F}_{\mathrm{MS}}$ and $\hat{g} = (z^2 - \hat{\sigma}^2)/z^2$. Expression (2.13) is Stein's (1981) unbiased risk estimator combined with an estimator of $\sigma^2$. Alternatively, the risk estimator $\hat{r}_{\mathrm{MS}}(D, f)$ follows from the argument for Mallows' (1973) $C_p$ criterion.

For fixed penalty matrix $D$, the *shrinkage-adaptive* MS($D$) estimator is defined to be $\hat{\eta}_{\mathrm{MS}}(D, \hat{f}_{\mathrm{MS}})$, where

$$\hat{f}_{\mathrm{MS}} = \underset{f \in \mathcal{F}_{\mathrm{MS}}}{\arg\min}\, \hat{r}_{\mathrm{MS}}(D, f) = \underset{f \in \mathcal{F}_{\mathrm{MS}}}{\arg\min}\, \mathrm{ave}[(f - \hat{g})^2 z]. \tag{2.14}$$

To accomplish the minimization, let $\mathcal{K} = \{k \in R^p : k_1 \geq k_2 \geq \cdots \geq k_p\}$ and let

$$\hat{k} = \underset{k \in \mathcal{K}}{\arg\min}\, \mathrm{ave}[(k - \hat{g})^2 z]. \tag{2.15}$$

Computation of $\hat{k}$ is a weighted isotonic least squares problem that can be solved in a finite number of steps with the pool-adjacent-violators algorithm [cf. Robertson, Wright and Dykstra (1988)]. Each component of $\hat{f}_{\mathrm{MS}}$ is then the positive part of the corresponding component of $\hat{k}$, as shown in Beran and Dümbgen (1998).

REMARK. The shrinkage adaptive PLS($D$) estimator is obtained by restricting the minimization in (2.14) to monotone shrinkage vectors of the form $f = f(\nu)$. This weighted nonlinear least squares computation is harder than constructing the more ambitious shrinkage adaptive MS($D$) estimator.

Useful in risk estimation is the high component variance estimator $\hat{\sigma}_{\mathrm{H}}^2$, which uses the strategy of pooling sums of squares from analysis of variance.

Choose $\overline{U}$ so that the concatenated matrix $(U|\overline{U})$ is orthogonal. Set $\bar{z} = \overline{U}'y$ in analogy to the earlier $z = U'y$. Then

$$(2.16) \quad \hat{\sigma}_H^2 = (n-q)^{-1}\left[\sum_{i=q+1}^{p} z_i^2 + |\bar{z}|^2\right] = (n-q)^{-1}\left[\sum_{i=q+1}^{p} z_i^2 + |y - \hat{\eta}_{LS}|^2\right],$$

where $q \leq \min\{p, n-1\}$. The bias of $\hat{\sigma}_H^2$ is $(n-q)^{-1}\sum_{i=q+1}^{p}\xi_i^2$. Consistency of $\hat{\sigma}_H^2$ is assured when this bias tends to zero as $n - q$ tends to infinity. When $q = p < n$, the estimator $\hat{\sigma}_H^2$ reduces to the least squares estimator $\hat{\sigma}_{LS}^2 = (n-p)^{-1}|y - \hat{\eta}_{LS}|^2$, which is unbiased. When $p = n$, the estimator $\hat{\sigma}^2$ is a pure pooling estimator whose bias is small if $(p-q)^{-1}\sum_{i=q+1}^{p}\xi_i^2$ is nearly zero. We will seek to arrange this through choice of the penalty matrix $D$.

2.3. *Adaptive ST estimators.* For $t \geq 0$ and $1 \leq i \leq p$, let $h_i(t, z) = [1 - t/|z_i|]_+$. Let

$$(2.17) \quad \mathcal{F}_{ST} = \mathcal{F}_{ST}(p) = \{f \in [0,1]^p : f_i = h_i(t, z) \text{ for } t \geq 0 \text{ and } 1 \leq i \leq p\}.$$

Unlike the monotone class $\mathcal{F}_{MS}$ defined in (2.8), the class $\mathcal{F}_{ST}$ of shrinkage vectors is data dependent. Let

$$(2.18) \quad \hat{\xi}_{ST}(D, f) = \{fz : f \in \mathcal{F}_{ST}\},$$

multiplication being performed componentwise as in S. The algebraic identity $\hat{h}_i(t, z)z_i = \text{sgn}(z_i)[|z_i| - t]_+$ connects $\hat{\xi}_{ST}(D, f)$ with the definition of soft-thresholding in Donoho and Johnstone (1995). The candidate ST estimators for $\eta$ associated with penalty matrix $D$ are

$$(2.19) \quad \hat{\eta}_{ST}(D, f) = U\hat{\xi}_{ST}(D, f) = U\,\text{diag}\{f\}U'y, \qquad f \in \mathcal{F}_{ST}.$$

Let $\widehat{G}$ denote the empirical *cumulative distribution function of the* $\{|z_i| : 1 \leq i \leq p\}$, let $G = \text{E}(\widehat{G})$ and define

$$(2.20) \quad r_{ST}(f, \xi, \sigma^2) = \sigma^2[1 - 2G(t)] + \int_0^\infty (u^2 \wedge t^2)\,dG(u), \qquad f \in \mathcal{F}_{ST},$$

where $\wedge$ denotes the minimum operator. It follows from Stein (1981) that the normalized quadratic risk of the candidate ST estimator is

$$(2.21) \quad R(\hat{\eta}_{ST}(D, f), \eta, \sigma^2) = r_{ST}(f, \xi, \sigma^2), \qquad f \in \mathcal{F}_{ST}.$$

Having devised a variance estimator $\hat{\sigma}^2$, we may estimate this risk by

$$(2.22) \quad \hat{r}_{ST}(D, f) = \hat{\sigma}^2[1 - 2\widehat{G}(t)] + \int_0^\infty (u \wedge t)^2\,d\widehat{G}(u), \qquad f \in \mathcal{F}_{ST}.$$

Let $t_p = (2\log(p))^{1/2}$. For fixed penalty matrix $D$, the *shrinkage-adaptive* $\mathrm{ST}(D)$ estimator is defined to be $\hat{\eta}_{\mathrm{ST}}(D, \hat{f}_{\mathrm{ST}})$, where

$$\hat{f}_{\mathrm{ST}} = h(\hat{t}, z) \qquad \text{where } \hat{t} = \underset{t \in [0, t_p]}{\arg\min}\, \hat{r}_{\mathrm{ST}}(D, t), \tag{2.23}$$

as in Donoho and Johnstone (1995). Because $\hat{t}$ must be one of the values $\{|z_i| : 1 \leq i \leq p\}$, it can be computed readily.

2.4. *Adaptive HS estimators.* Let $p_1 = \lfloor \alpha p \rfloor$, where $\lfloor \cdot \rfloor$ denotes integer part and the *split* fraction $\alpha \in [0, 1]$. For any vector $k \in R^p$, define the subvectors $k_{(1)} = \{k_i : 1 \leq i \leq p_1\}$ and $k_{(2)} = \{k_i : p_1 + 1 \leq i \leq p\}$ of respective dimensions $p_1$ and $p_2 = p - p_1$. Candidate HS estimators apply separate shrinkage strategies to the subvectors $z_{(1)}$ and $z_{(2)}$ of $z$. We focus on the $\mathrm{MS} \times \mathrm{ST}$ hybrid because it proves particularly effective in the examples to be considered. The definitions of the $\mathrm{MS} \times \mathrm{MS}$, $\mathrm{ST} \times \mathrm{ST}$ and of $\mathrm{ST} \times \mathrm{MS}$ hybrids are analogous.

Efromovich (1999) considered HS of wavelet coefficients in which MS is replaced by a certain linear shrinkage methodology and ST is replaced by hard-thresholding. In both that paper and here, the aim is to compromise beneficially between a shrinkage approach that assumes regression coefficients are ordered in importance and a shrinkage approach that relies on sparsity of important regression coefficients. Considerable technical differences exist. We apply adaptive MS rather than Efromovich–Pinsker shrinkage to the low-frequency regression coefficients. On the remaining coefficients, we use ST rather than hard-thresholding and select the soft-threshold to minimize estimated risk. The regularity conditions for Stein's (1981) risk estimator are satisfied by soft-thresholding but not by hard-thresholding.

Let

$$\mathcal{F}_{\mathrm{HS}} = \{f : f_{(1)} \in \mathcal{F}_{\mathrm{MS}}(p_1), f_{(2)} \in \mathcal{F}_{\mathrm{ST}}(p_2)\} \tag{2.24}$$

and let

$$\hat{\xi}_{\mathrm{HS}}(D, \alpha, f) = fz, \qquad f \in \mathcal{F}_{\mathrm{HS}}. \tag{2.25}$$

The candidate $\mathrm{MS} \times \mathrm{ST}$ HS estimators for $\eta$ associated with penalty matrix $D$ are defined by

$$\hat{\eta}_{\mathrm{HS}}(D, \alpha, f) = U\hat{\xi}_{\mathrm{HS}}(D, \alpha, f) = U\,\mathrm{diag}\{f\}U'y, \qquad f \in \mathcal{F}_{\mathrm{HS}}. \tag{2.26}$$

From the preceding sections, it follows that the normalized quadratic risk of this candidate HS estimator is

$$\begin{aligned}
R(\hat{\eta}_{\mathrm{HS}}(D, \alpha, f), \eta, \sigma^2) \\
= p^{-1}[p_1 r_{\mathrm{MS}}(f_{(1)}, \xi_{(1)}, \sigma^2) + p_2 r_{\mathrm{ST}}(f_{(2)}, \xi_{(2)}, \sigma^2)], \qquad f \in \mathcal{F}_{\mathrm{HS}}.
\end{aligned} \tag{2.27}$$

TABLE 1

| LS | PLS(D4) | MS(D4) | ST(D4) | HS(D4) |
|---|---|---|---|---|
| 0.0115 | 0.0093 | 0.0071 | 0.0047 | 0.0039 |

Write $\hat{r}_{\mathrm{MS}}(D, f_{(1)})$ for the risk estimator (2.13) computed on the subvector $z_{(1)}$. Similarly, write $\hat{r}_{\mathrm{ST}}(D, f_{(2)})$ for the risk estimator (2.22) computed on the subvector $z_{(2)}$. The risk of the candidate HS estimator is then estimated by

$$(2.28) \quad \hat{r}_{\mathrm{HS}}(D, \alpha, f) = p^{-1}[p_1 \hat{r}_{\mathrm{MS}}(D, f_{(1)}) + p_2 \hat{r}_{\mathrm{ST}}(D, f_{(2)})], \qquad f \in \mathcal{F}_{\mathrm{HS}}.$$

For fixed penalty matrix $D$ and split fraction $\alpha$, the *shrinkage-adaptive* HS($D$) estimator is defined to be $\hat{\eta}_{\mathrm{HS}}(D, \alpha, \hat{f}_{\mathrm{HS}})$, where

$$(2.29) \qquad \hat{f}_{\mathrm{HS}} = \underset{f \in \mathcal{F}_{\mathrm{HS}}}{\arg \min} \, \hat{r}_{\mathrm{HS}}(D, \alpha, f).$$

The minimization is accomplished by minimizing separately each of the two summands on the right-hand side of (2.28) in the manner discussed previously.

2.5. *A case study.* Figure 1 presents competing fits to monthly Australian red wine sales (in kiloliters) from January 1980 to October 1991. The data are taken from Brockwell and Davis (1996) and the ordinal factor is month. Here $n = p = 142$. The penalty matrix is the fourth difference operator D4, which is defined explicitly in Section 3.2. The high component variance estimate $\hat{\sigma}_{\mathrm{H}}^2$ is determined by (2.16) with $q = \lfloor 0.85p \rfloor$. The partition in the definition of HS(D4) uses $\alpha = 0.3$. Adaptation to minimize estimated risk selected the values of $\alpha$ and of the penalty matrix from a class of possibilities described in Section 3.3. The estimated risks of the competing estimators are shown in Table 1.

The LS fit (not shown) coincides with the raw data. On the basis of estimated risk, PLS(D4) is only a modest improvement over LS, MS(D4) is preferable, while ST(D4) and HS(D4) are substantially preferable, the hybrid estimator being best. Theorem 4.1 shows that, under model (1.1), the estimated risks of these adaptive estimators approximate their risks under the saturated model as $p$ tends to infinity.

On looking closely at Figure 1, we discern a regular seasonal pattern in the HS(D4) and ST(D4) fits. Each year, estimated mean monthly red wine sales rise steadily from an annual low in January to a peak around July or August (winter in Australia) and then drop into a trough with a secondary peak around November or December (in time for the Christmas holiday season). The adaptive fits with smallest estimated risk have recovered a highly intelligible seasonal pattern in sales that may be linked to seasonal patterns in market demand and in winery operations after harvest and fermentation.

Figure 2 examines what is going on behind the fits. The residuals from the HS(D4) fit are plausibly homoscedastic. A Q–Q plot (not shown) indicates that their marginal distribution is roughly normal, apart from outliers. This illustrates the tendency of our procedures to fit the data in terms of the motivating model. Subplot $(1,2)$ plots the transformed components $\{|z_i|^{1/2}\operatorname{sgn}(z_i) : 1 \le i \le p\}$ of the coefficients $z = U'y$. The square root transformation reduces the vertical range of the plot and makes more visible the behavior of small components of $z$. Evidently, the first four columns of $U$ are crucial in representing $y$ and so $\eta$. Blips in this plot at certain higher-order components suggest that the corresponding basis vectors may also be important in estimating $\eta$. We call subplot $(1,2)$ an *empirical basis economy plot*. The concept of basis economy is treated formally in Section 3.1. As well, this subplot suggests the choice of $q$ that enters into the high-component variant estimator $\hat{\sigma}_H^2$.

The four shrinkage vector subplots in Figure 2 display the shrinkage vectors that define the competing adaptive fits. Because the shrinkage vectors of the PLS(D4) and MS(D4) estimates are necessarily monotone, both give considerable weight to many components of $z$ so as not to disregard the small blips discussed above. The ST(D4) and HS(D4) estimates are better able to select the more important components of $z$, thereby reducing estimated risk through tradeoff of estimated variance against bias. Note that the HS(D4) estimate disregards more of the higher-order components of $z$ than does ST(D4).

2.6. *Experiments with artificial data.* Figures 3 and 4 exhibit the competing adaptive estimators on two sets of artificial monthly data that bracket the situation found in Example 1. In this experiment, $p = n = 200$, the factor levels are $\{s_i = i : 1 \le i \le 200\}$, and the means at which we have one noisy observation are

Smooth: $m_1(s_i) = 2 - 50((s_i/200 - 0.25)(s_i/200 - 0.75))^2$,
Very Wiggly: $m_2(s_i) = m_1(s_i/200) - 0.25\sin(100\pi(s_i/200))$.

The observations are given by $y_i = m(s_i) + e_i$, where the $\{e_i\}$ form a single pseudo-random sample drawn from the $N(0, \sigma^2 I_{200})$ distribution with $\sigma = 0.5$. In the data analysis, the variance $\sigma^2$ is estimated by the high component estimator $\hat{\sigma}_H^2$ defined in (2.16), with $q = 0.75p$.

Fitting this artificial data is a one-way layout problem rather than a curve estimation problem because the measurements are deemed to be monthly as in Example 1. The means in the Smooth case vary more slowly than those estimated in Example 1, while the means in the Very Wiggly case vary more rapidly. The goal is to learn how the competing adaptive estimators perform

HYBRID SHRINKAGE IN ONE-WAY ANOVA 17in both scenarios. The first rows of Figures 3 and 4 give the scatterplots of the Smooth and Very Wiggly data, respectively. To the human eye, these scatterplots are scarcely distinguishable. Good estimators of the unknown mean vectors seek to do better than the eye.

The penalty matrix used for both sets of artificial data is the fourth difference operator D4. The partition in the definition of HS(D4) uses $\alpha = 0.05$. Adaptation to minimize estimated risk selected these values of $\alpha$ and of the penalty matrix from a class of possibilities described in Section 3.3. According to the asymptotics in Section 4, the risk, loss and estimated risk all converge to a common limit. In the present experiment with artificial data, the losses are readily computed. For the Smooth data, the estimated risks and actual losses of the competing estimators are shown in Table 2.

We note that the estimated risks for the shrinkage adaptive estimators are negative. The actual losses are small and convergence to asymptotic limits has not happened. Nevertheless, the estimated risks reflect the ordering of the true losses. In Figure 3 the visual quality of the competing fits follows the same ordering. The interpolated ST(D4) estimate is unsatisfactorily jagged, though certainly better than the LS estimate. The MS(D4) and HS(D4) estimates are close to the truth, though the latter exhibits a small ripple not present in the actual mean vector. The basis economy plot in the last subplot of Figure 3 suggests that the D4 penalty basis is economical in this example. This is verified by examining the corresponding plot of $\xi$ (not shown) computed from the true mean function.

For the Very Wiggly data, the estimated risks and actual losses of the competing estimators are shown in Table 3.

In Figure 4, interpolating lines have been added to guide the eye through the sequence of means or estimated means. They have no further significance because we are not doing curve estimation. The HS(D4) estimate is best visually, as well as in loss. The HS(D4) and ST(D4) estimates both indicate the amplitude of the high frequency component in the unknown mean more successfully than the MS(D4) estimate. However, the actual loss of the ST(D4) estimate exceeds that of the MS(D4) estimate. Both casual scrutiny and the ordering of the estimated losses make ST(D4) look better than it is. Evidently the asymptotics have not fully taken hold. The basis economy plot in subplot $(3, 2)$ of Figure 4 reveals the possible importance of component $z_{102}$. In the Very Wiggly case, the D4 penalty basis is sparse in the sense that most components of $\xi$ are small. However, it is not economical because a high-order basis vector is needed to approximate the high frequency sinusoidal component in the mean.

Table 2

|  | LS | MS(D4) | ST(D4) | HS(D4) |
|---|---|---|---|---|
| Estimated risk | 0.2846 | −0.0434 | −0.0296 | −0.0449 |
| Loss | 0.2325 | 0.0072 | 0.0358 | 0.0077 |

Table 3

|  | LS | MS(D4) | ST(D4) | HS(D4) |
|---|---|---|---|---|
| Estimated risk | 0.2842 | −0.0063 | −0.0239 | −0.0350 |
| Loss | 0.2325 | 0.0313 | 0.0447 | 0.0285 |

In this experiment, the HS(D4) estimate, unlike the others considered, performs well in *both* the Smooth case and the Very Wiggly case. This is empirical evidence in its favor.

**3. Penalty matrix and split fraction.** For monotone shrinkage, the ideal choice of basis $U$ would have its first column proportional to the unknown mean vector $\eta$ so that only the first component of $\xi = U'\eta$ is nonzero. Then the choice of shrinkage vector $f$ to minimize risk would have first component equal to 1 and all other components equal to 0. Though unrealizable, this ideal choice indicates that prior information or conjecture about $\eta$ should be exploited in selecting $U$. We say informally that the columns of $U$ provide an *economical* basis for the regression space if all but the first few components of $\xi$ are very nearly zero. Construction of the basis $U$ via a penalty matrix $D$—the method used in this paper—is a practical way of using vague prior information or conjecture about the function $m$ to find a plausibly economical basis for expressing the mean vector $\eta$.

3.1. *The role of basis economy.* Mathematical analysis of an idealized economy concept reveals the importance of basis economy in reducing risk through monotone shrinkage. For every $b \in (0, 1]$, let $\mathcal{E}_M(b) = \{a \in R^p : a_i = 1 \text{ if } 1 \leq i \leq bp, 1 \leq a_{\lfloor bp \rfloor + 1} \leq \cdots \leq a_p \leq \infty\}$. For every $a \in \mathcal{E}_M(b)$, every $r > 0$ and every $\sigma^2 > 0$, define the ellipsoid

$$(3.1) \qquad E(r, a, \sigma^2) = \{\xi \in R^p : \operatorname{ave}(a\xi^2) \leq \sigma^2 r\}.$$

If $\xi \in E(r, a, \sigma^2)$ and $a_i = \infty$, it is to be understood that $\xi_i = 0$ and $a_i^{-1} = 0$. We consider bases $U$ such that, in the resulting canonical model, $\xi \in E(r, a, \sigma^2)$ for some $r > 0$, some $a \in \mathcal{E}_M(b)$ and some $b \in (0, 1]$.

A finite-dimensional specialization of Pinsker's (1980) theorem, given by Beran and Dümbgen (1998), implies the next theorem on asymptotic minimaxity of adaptive MS estimators of $\eta$. The proof follows from the discussion in Section 4 of Beran (2000). Let $\xi_0^2 = \sigma^2[(\gamma/a)^{1/2} - 1]_+$, where $\gamma$ is the unique positive number such that $\mathrm{ave}(\xi_0^2) = \sigma^2 r$. Define

$$\nu_p(r, a, \sigma^2) = \sigma^2 \, \mathrm{ave}[\xi_0^2 / (\sigma^2 + \xi_0^2)]. \tag{3.2}$$

THEOREM 3.1. *Fix the penalty basis $U$ by choice of $D$ or otherwise. For every $b \in (0, 1]$, every $a \in \mathcal{E}_{\mathrm{M}}(b)$, every $r > 0$ and every $\sigma^2 > 0$,*

$$\lim_{p \to \infty} \left[ \inf_{\hat{\eta}} \sup_{\xi \in E(r,a,\sigma^2)} R(\hat{\eta}, \eta, \sigma^2) / \nu_p(r, a, \sigma^2) \right] = 1. \tag{3.3}$$

*The shrinkage-adaptive estimator $\hat{\eta}_{\mathrm{MS}}(D, \hat{f}_{\mathrm{MS}})$ achieves asymptotic minimax bound* (3.3) *in that*

$$\lim_{p \to \infty} \left[ \sup_{\xi \in E(r,a,\sigma^2)} R(\hat{\eta}_{\mathrm{MS}}(D, \hat{f}_{\mathrm{MS}}), \eta, \sigma^2) / \nu_p(r, a, \sigma^2) \right] = 1. \tag{3.4}$$

What does this theorem tell us? First, note that the asymptotic minimax risk $\nu_p(r, a, \sigma^2)$ in (3.3) is monotone decreasing in the vector $a$. Thus, if $\xi = U'\eta \in E(r, a, \sigma^2)$ for relatively small $b$ and relatively large vector $a$—in other words, if the basis is economical for expressing $\eta$—then the asymptotic minimax risk is relatively small compared to the risk $\sigma^2$ of the LS estimator. Second, (3.4) indicates that the adaptive MS estimator achieves the asymptotic minimax risk for every degree of basis economy. Even a poor choice of basis for adaptive MS estimation does not lead to disaster relative to LS estimation.

A special case of Theorem 3.1 makes both points obvious, albeit in a simplified setting. Let $\mathcal{B}(b) = \{a \in \mathcal{E}_{\mathrm{M}}(b) : a_i = \infty \text{ if } \lfloor bp \rfloor + 1 \leq i \leq p\}$. In Theorem 3.1, replacing $a \in \mathcal{E}_{\mathrm{M}}(b)$ with the stronger restriction $a \in \mathcal{B}(b)$ and $\nu_p(r, a, \sigma^2)$ with the evaluation $\sigma^2 rb/(r + b)$ gives a valid statement. In this simplified setting, basis economy corresponds to a small value of $b$. The ratio of the asymptotic minimax risk to the risk of the LS estimator is small whenever $b$ is small; and the adaptive MS estimator is still asymptotically minimax.

3.2. *Local annihilators.* Difference operators are well-established as penalty matrices for PLS when the ordinal factor levels $s = (s_1, s_2, \ldots, s_p)$, with $s_1 < s_2 < \cdots < s_p$, are equally spaced [cf. Press, Teukolsky, Vetterling and Flannery (1992), Section 18.5]. To define the $d$th difference matrix $D_d$, first define the $(p-1) \times p$ matrix $\Delta(p) = \{\delta_{i,j}\}$, in which $\delta_{i,i} = 1$, $\delta_{i,i+1} = -1$ for every $i$ and all other entries are zero. Then

$$D_1 = \Delta(p), \qquad D_d = \Delta(p - d + 1) D_{d-1} \qquad \text{for } 2 \leq d < p. \tag{3.5}$$

Evidently, the $(p-d) \times p$ matrix $D_d$ annihilates powers of $s$ up to power $d-1$ in the sense that

$$(3.6) \qquad D_d s^k = 0 \qquad \text{for } 0 \leq k \leq d-1.$$

Here $s^k$ denotes the column vector $(s_1^k, \ldots, s_p^k)'$. Moreover, in row $i$ of $D_d$, the elements not in columns $i, i+1, \ldots, i+d$ are zero.

Suppose for simplicity that $X = I_p$. Let $U$ be the penalty basis generated by penalty matrix $D_d$. By the variational characterization of $U$ given in Section 2.1, the space spanned by the first $d$ columns of $U$ consists of vectors $v$ that satisfy $D_d v = 0$. When $m$ behaves locally like a polynomial of degree $d-1$ and the value of $d$ is modest, then this penalty basis is economical for $\eta$. Such considerations support the use of difference operators as candidate penalty matrices when the factor levels are equally spaced.

When $m$ is expected to behave locally like a polynomial of degree $d-1$, but the factor levels in $s$ are not equally spaced, we replace $D_d$ as follows. For every integer $1 \leq d < p$, the *local polynomial annihilator* $A_d$ is a $(p-d) \times p$ matrix characterized through three conditions. First, for every possible $i$, all elements in the $i$th row of $A_d$ other than $\{a_{i,j} : i \leq j \leq i+d\}$ are zero. Second, $A_d$ satisfies the orthogonality conditions

$$(3.7) \qquad A_d s^k = 0 \qquad \text{for } 0 \leq k \leq d-1.$$

Third, each row vector in $A_d$ has unit length. These requirements are met by setting the nonzero elements in the $i$th row of $A_d$ equal to the basis vector of degree $d$ in the orthonormal polynomial basis that is defined on the $d+1$ design points $(s_i, \ldots, s_{i+d})$. The S-Plus function `poly` accomplishes this computation. When the components of $s$ are equally spaced, $A_d$ is just a scalar multiple of the $d$th difference matrix $D_d$.

3.3. *Adaptive choice of penalty matrix and split.* As we have seen, a penalty basis ideally exploits, through choice of the penalty matrix, informed conjecture about the function $m$ in (1.1). When this is the case, penalty bases are often reasonably economical. However, if the prior information is weak or flawed, some of the higher-order components of $\xi$ may not be negligible. Soft-thresholding handles possibly isolated higher-order components of $\xi$ that need to be considered in the fit. The choice of dividing point $p_1$ between monotone shrinkage and soft-thresholding in the MS × ST HS estimator then becomes important. We will use the strategy of minimizing estimated risk to select $D$ and $p_1$, in addition to the shrinkage vectors.

Given a set $\mathcal{D}$ of candidate penalty matrices, such as $\{A_d : 1 \leq d \leq k\}$, we select an empirically best MS estimator as follows. Over shrinkage class $\mathcal{F}_{\mathrm{MS}}$ and over penalty matrix class $\mathcal{D}$, the fully adaptive MS estimator of $\eta$ is defined to be $\hat{\eta}_{\mathcal{D},\mathrm{MS}} = \hat{\eta}_{\mathrm{MS}}(\widehat{D}, \hat{f})$, where

$$(3.8) \qquad (\widehat{D}, \hat{f}) = \underset{D \in \mathcal{D}, f \in \mathcal{F}_{\mathrm{MS}}}{\arg\min} \; \hat{r}(D, f).$$

The fully adaptive ST estimator $\hat{\eta}_{\mathcal{D},\text{ST}}$ is defined analogously, replacing $\mathcal{F}_{\text{MS}}$ in (3.8) with $\mathcal{F}_{\text{ST}}$.

For HS estimators, it is also desirable to explore competing choices of $p_1 = \lfloor \alpha p \rfloor$, where $\lfloor \cdot \rfloor$ denotes integer part and candidate values of $\alpha$ lie in a specified subset $\mathcal{A}$ of $[0,1]$. Over shrinkage class $\mathcal{F}_{\text{HS}}$, over penalty matrix class $\mathcal{D}$ and over split fraction class $\mathcal{A}$, the fully adaptive HS estimator of $\eta$ is defined to be $\hat{\eta}_{\mathcal{D},\mathcal{A},\text{HS}} = \hat{\eta}_{\text{HS}}(\hat{D}, \hat{\alpha}, \hat{f})$, where

$$(3.9) \qquad (\hat{D}, \hat{\alpha}, \hat{f}) = \operatorname*{arg\,min}_{D \in \mathcal{D}, \alpha \in \mathcal{A}, f \in \mathcal{F}_{\text{HS}}} \hat{r}(D, \alpha, f).$$

The asymptotics in Section 4 support choosing $\alpha$ and $D$ to minimize estimated risk provided the cardinalities of $\mathcal{A}$ and of $\mathcal{D}$ grow slowly as $p$ increases. The numerical examples in Sections 2.5 and 2.6 used $\mathcal{D} = \{D_d : 1 \leq d \leq 6\}$ and $\mathcal{A} = \{0.05k : 0 \leq k \leq 20\}$. The asymptotics given do not care whether the candidate bases are constructed as penalty bases. However, minimizing estimated risk over a very large class of bases should not be expected to yield a good estimator of $\eta$. For instance, the MS estimator that minimizes the estimated risk of $U \operatorname{diag}\{f\} U' y$ over all $f \in \mathcal{F}_{\text{MS}}$ and over all permutations of the columns of a fixed basis matrix $U$ is dominated by the LS estimator in the saturated model. Remark A on page 1829 of Beran and Dümbgen (1998) gives a proof. In such cases, the covering numbers used in the asymptotics of Section 4 are too large for Theorem 4.1 to hold.

**4. Asymptotics of adaptation.** The main purpose of this section is to analyze the asymptotic loss and risk of the adaptive $\text{ST}(D)$ and $\text{HS}(D)$ estimators under the saturated Gaussian one-way layout. The results build on techniques developed by Beran and Dümgben (1998) for adaptive $\text{MS}(D)$ estimators. First we show that minimizing estimated risk over shrinkage class $\mathcal{F}_{\text{MS}}$ or $\mathcal{F}_{\text{ST}}$ for fixed penalty matrix $D$ succeeds in minimizing risk asymptotically over that shrinkage class as the dimension $p$ of the regression space tends to infinity. Moreover, the estimated risk of the adaptive estimator converges to its actual loss and risk. In this fashion, estimated risks provide a credible tool for ranking competing shrinkage estimators. Second, we provide conditions under which simultaneous adaptation over shrinkage class $\mathcal{F}_{\text{HS}}$, over penalty matrix class $\mathcal{D}$ and over split fraction class $\mathcal{A}$ works in the senses just described. The results require no smoothness assumptions on the unknown mean vector $\eta$.

4.1. *Adaptation works.* For any vector $h \in R^p$, let $\|h\| = \max_{1 \leq i \leq p} |h_i|$. The generic subscript $\mathcal{F}$ stands for $\mathcal{F}_{\text{MS}}$ or $\mathcal{F}_{\text{ST}}$ or $\mathcal{F}_{\text{HS}}$, according to the choice of candidate estimator class.

THEOREM 4.1. *Let $\mathcal{F}$ be either $\mathcal{F}_{\mathrm{MS}}$ or $\mathcal{F}_{\mathrm{ST}}$. Suppose that $\hat{\sigma}^2$ is consistent in that, for every $c > 0$ and $\sigma^2 > 0$,*

$$(4.1) \qquad \lim_{p \to \infty} \sup_{\|\xi\| \leq c} \mathrm{E}|\hat{\sigma}^2 - \sigma^2| = 0.$$

(a) *Let $V(f)$ denote either the loss $L(\hat{\eta}_{\mathcal{F}}(D,f), \eta)$ or the estimated risk $\hat{r}_{\mathcal{F}}(D,f)$. Then, for every penalty matrix $D$, every $c > 0$ and every $\sigma^2 > 0$,*

$$(4.2) \qquad \lim_{p \to \infty} \sup_{\|\xi\| \leq c} \mathrm{E} \sup_{f \in \mathcal{F}} |V(f) - R(\hat{\eta}_{\mathcal{F}}(D,f), \eta, \sigma^2)| = 0.$$

(b) *If $\hat{f} = \arg\min_{f \in \mathcal{F}} \hat{r}_{\mathcal{F}}(D,f)$, then*

$$(4.3) \qquad \lim_{p \to \infty} \sup_{\|\xi\| \leq c} \left| R(\hat{\eta}_{\mathcal{F}}(D,\hat{f}), \eta, \sigma^2) - \min_{f \in \mathcal{F}} R(\hat{\eta}_{\mathcal{F}}(D,f), \eta, \sigma^2) \right| = 0.$$

(c) *For $W$ equal to either $L(\hat{\eta}_{\mathcal{F}}(D,\hat{f}), \eta)$ or $R(\hat{\eta}_{\mathcal{F}}(D,\hat{f}), \eta, \sigma^2)$,*

$$(4.4) \qquad \lim_{p \to \infty} \sup_{\|\xi\| \leq c} \mathrm{E}|\hat{r}_{\mathcal{F}}(D,\hat{f}) - W| = 0.$$

(d) *Let $\#\mathcal{D}$ denote the cardinality of $\mathcal{D}$. Convergences (4.2) to (4.4) hold for the fully adaptive MS estimator $\hat{\eta}_{\mathcal{D},\mathrm{MS}}$ defined through (3.8) if $(\#\mathcal{D})p^{-1/2}$ and $(\#\mathcal{D})\mathrm{E}|\hat{\sigma}^2 - \sigma^2|$ both tend to zero as $p \to \infty$. They hold for the fully adaptive ST estimator if $(\#\mathcal{D})p^{-1/2}(\log(p))^{1/4}$ and $(\#\mathcal{D})\mathrm{E}|\hat{\sigma}^2 - \sigma^2|$ both tend to zero as $p \to \infty$.*

(e) *Convergences (4.2) to (4.4) hold for the fully adaptive HS estimator $\hat{\eta}_{\mathcal{D},\mathcal{A},\mathrm{HS}}$ defined in (3.8) if $\max\{\#\mathcal{A}, \#\mathcal{D}\}p^{-1/2}(\log(p))^{1/4}$ and $\max\{\#\mathcal{A}, \#\mathcal{D}\}\mathrm{E}|\hat{\sigma}^2 - \sigma^2|$ both tend to zero as $p \to \infty$.*

Parts (a)–(c) refer to the case of fixed $D$. By part (a), the loss, risk and estimated risk of a candidate estimator converge together, uniformly over $\mathcal{F} = \mathcal{F}_{\mathrm{MS}}$ or $\mathcal{F}_{\mathrm{ST}}$. This makes the estimated risk of candidate estimators indexed by $\mathcal{F}$ a trustworthy surrogate for true risk or loss. By part (b), the risk of the shrinkage-adaptive estimator $\hat{\eta}_{\mathcal{F}}(D, \hat{f})$ converges to that of the best candidate estimator. Part (c) shows that the loss, risk and plug-in estimated risk of an adaptive estimator converge together asymptotically. Part (d) extends these findings to MS and ST estimators that adapt over both $f$ and $D$. Part (e) does the same for HS estimators that adapt over $f$, $D$ and $\alpha$.

Condition (4.1) holds for the variance estimator $\hat{\sigma}^2_{\mathrm{LS}}$ if $n - p$ tends to infinity with $p$. Asymptotic results for other variance estimators are given in Beran (1996) and Beran and Dümbgen (1998).

4.2. *Auxiliary result.* The proof of Theorem 4.1 uses techniques from empirical process theory. Theorem 4.2 below is taken from Beran and Dümbgen (1998). It follows from standard symmetrization arguments and Pisier's (1983) form of the chaining lemma [see also Pollard (1990), Sections 2 and 3]. Let $S = \sum_{i=1}^{p} \phi_i$, where $\phi_1, \phi_2, \ldots, \phi_p$ are independent stochastic processes on an index set $\mathcal{T}$. All $\phi_i$ have continuous sample paths with respect to some metric $\rho$ on $\mathcal{T}$ such that $(\mathcal{T}, \rho)$ is separable. Define a random pseudo-metric $\hat{m}$ on $\mathcal{T}$ through

$$(4.5) \qquad \hat{m}^2(s,t) = \sum_{i=1}^{p}[\phi_i(s) - \phi_i(t)]^2.$$

For any pseudo-metric $\nu$ on $\mathcal{T}$, define the covering numbers

$$(4.6) \qquad N(u, \mathcal{T}, \nu) = \min\left\{\#\mathcal{T}_0 : \mathcal{T}_0 \subset \mathcal{T}, \inf_{t_0 \in \mathcal{T}_0} \nu(t_0, t) \leq u \ \forall t \in \mathcal{T}\right\}.$$

THEOREM 4.2. *Suppose that $S(t_1) \equiv 0$ for some $t_1 \in \mathcal{T}$. Then there exists a finite constant $C > 0$ such that*

$$(4.7) \qquad \mathrm{E}\sup_{t \in \mathcal{T}}|S(t) - \mathrm{E}S(t)| \leq C\mathrm{E}\int_0^{\widehat{D}} \log^{1/2}[N(u, \mathcal{T}, \hat{m})]\,du,$$

*where $\widehat{D} = \sup_{t \in \mathcal{T}} \hat{m}(t, t_1)$.*

4.3. *Proof of Theorem* 4.1. The portion of Theorem 4.1 that concerns $\mathcal{F} = \mathcal{F}_{\mathrm{MS}}$ follows from results in Section 6 of Beran and Dümbgen (1998). We continue by proving parts (a)–(c) for $\mathcal{F} = \mathcal{F}_{\mathrm{ST}}$. For this discussion of soft-thresholding, let $\mathcal{T} = [0, t_p]$ with $t_p = (2\log(p))^{1/2}$.

(a) Suppose that $V(f) = \hat{r}_{\mathrm{ST}}(D, f)$ for $f \in \mathcal{F}_{\mathrm{ST}}$. In view of (2.22) and (4.1), it suffices to show that

$$(4.8) \qquad \lim_{p \to \infty} \sup_{\|\xi\| \leq c} \mathrm{E}\sup_{t \in \mathcal{T}} |\widehat{G}(t) - G(t)| = 0$$

and

$$(4.9) \qquad \lim_{p \to \infty} \sup_{\|\xi\| \leq c} \mathrm{E}\sup_{t \in \mathcal{T}} \left|\int_0^\infty (u^2 \wedge t^2)\,d[\widehat{G}(u) - G(u)]\right| = 0.$$

In Theorem 4.2, take $\phi_i(t) = p^{-1}I(|z_i| \leq t)$. Then $S(t) = \widehat{G}(t)$, $\hat{m}^2(s,t) = p^{-1}|\widehat{G}(t) - \widehat{G}(s)|$, $t_1 = 0$, $\widehat{D} = p^{-1/2}$, and

$$(4.10) \qquad \begin{aligned} N(u, \mathcal{T}, \hat{m}) &= \min\left\{\#\mathcal{T}_0 : \mathcal{T}_0 \subset \mathcal{T}, \inf_{t_0 \in \mathcal{T}_0} \hat{m}^2(t_0, t) \leq u^2 \ \forall t \in \mathcal{T}\right\} \\ &\leq 1 + (pu^2)^{-1}. \end{aligned}$$

Then

$$\int_0^{\widehat{D}} \log^{1/2}[N(u, \mathcal{T}, \hat{m})] \, du \leq \int_0^{p^{-1/2}} \log^{1/2}[1 + (pu^2)^{-1}] \, du$$

(4.11)

$$= p^{-1/2} \int_0^1 \log^{1/2}(1 + v^{-2}) \, dv.$$

Because the rightmost integral is finite, (4.11) and (4.7) imply

(4.12) $$\mathrm{E} \sup_{t \in \mathcal{T}} |\widehat{G}(t) - G(t)| \leq C p^{-1/2}.$$

Limit (4.8) follows.

Next, observe that

$$\int_0^\infty (u^2 \wedge t^2) \, d\widehat{G}(u) = p^{-1} \sum_{i=1}^p z_i^2 I(|z_i| \leq t) + p^{-1} \sum_{i=1}^p t^2 I(|z_i| > t)$$

(4.13)

$$= S_1(t) + S_2(t), \quad \text{say}.$$

To analyze $S_1(t)$, let $\phi_i(t) = p^{-1} z_i^2 I(|z_i| \leq t)$. For any integer $r \geq 1$, let

(4.14) $$a_r = p^{-1} \sum_{i=1}^p |z_i|^r.$$

Now, using Cauchy–Schwarz, $\hat{m}^2(s, t) \leq p^{-1} a_8^{1/2} |\widehat{G}(t) - \widehat{G}(s)|^{1/2}$ and $\widehat{D} \leq p^{-1/2} a_8^{1/4}$. By reasoning akin to that in (4.10),

(4.15) $$N(u, \mathcal{T}, \hat{m}) \leq 1 + a_8(p^2 u^4)^{-1}.$$

Consequently, by (4.7) and a calculation like that in (4.11),

$$\mathrm{E} \sup_{t \in \mathcal{T}} |S_1(t) - \mathrm{E} S_1(t)| \leq C p^{-1/2} \mathrm{E} a_8^{1/4} \int_0^1 \log^{1/2}(1 + v^{-4}) \, dv$$

(4.16)

$$\leq C' p^{-1/2} \mathrm{E} a_8^{1/4}.$$

To analyze $S_2(t)$, let $\phi_i(t) = p^{-1} t^2 I(|z_i| > t)$. If $s \leq t$,

$$p^2 [\phi_i(s) - \phi_i(t)]^2 = [(s^2 - t^2) I(|z_i| > t) + s^2 I(s < |z_i| \leq t)]^2$$

(4.17)
$$\leq 2(s^2 - t^2)^2 I(|z_i| > t) + 2 s^4 I(s < |z_i| \leq t)$$

$$\leq 8 z_i^2 (s - t)^2 I(|z_i| > t) + 2 z_i^4 I(s < |z_i| \leq t).$$

Similarly for $t \leq s$. From this and Cauchy–Schwarz, $\hat{m}^2(s, t) \leq \hat{m}_1^2(s, t) + \hat{m}_2^2(s, t)$, where

$$\hat{m}_1^2(s, t) = p^{-1} 8 (s - t)^2 a_4^{1/2} [1 - \widehat{G}(\max(s, t))]^{1/2},$$

(4.18)

$$\hat{m}_2^2(s, t) = p^{-1} 2 a_8^{1/2} |\widehat{G}(s) - \widehat{G}(t)|^{1/2}.$$

By the first line in (4.17), $\widehat{D} \leq Ap^{-1/2}a_4^{1/2} \leq Ap^{-1/2}a_8^{1/4}$ for some finite constant $A$. Moreover,

(4.19)
$$N(u, \mathcal{T}, \hat{m}_1) \leq 1 + 8^{1/2}a_4^{1/4}t_p(p^{1/2}u)^{-1},$$
$$N(u, \mathcal{T}, \hat{m}_2) \leq 1 + 4a_8(p^2u^4)^{-1}$$

by reasoning similar to that for $S_1(t)$.

Because $\hat{m}(s,t) \leq \hat{m}_1(s,t) + \hat{m}_2(s,t)$,

(4.20) $\quad N(u, \mathcal{T}, \hat{m}) \leq 2\max\{N(u/2, \mathcal{T}, \hat{m}_1), N(u/2, \mathcal{T}, \hat{m}_2)\}$

and so

(4.21)
$$\int_0^{\widehat{D}} \log^{1/2}[N(u, \mathcal{T}, \hat{m})]\, du$$
$$\leq 2^{1/2} \int_0^{\widehat{D}} \log^{1/2}[N(u/2, \mathcal{T}, \hat{m}_1)]\, du$$
$$+ 2^{1/2} \int_0^{\widehat{D}} \log^{1/2}[N(u/2, \mathcal{T}, \hat{m}_2)]\, du.$$

The expectation of the second integral on the right-hand side is bounded from above by a constant times $p^{-1/2}$, as in (4.16). The expectation of the first integral on the right-hand side is bounded from above by a constant times $p^{-1/2}t_p^{1/2}$. Hence, by Theorem 4.2,

(4.22) $\quad \mathrm{E}\sup_{t \in \mathcal{T}}|S_2(t) - \mathrm{E}S_2(t)| \leq C_1''p^{-1/2} + C_2''p^{-1/2}t_p^{1/2}.$

Limit (4.9) now follows from (4.16) and (4.22). This establishes (4.2) for $V(f) = \hat{r}_{\mathrm{ST}}(D, f)$.

Next, suppose that $V(f) = L(\hat{\eta}_{\mathrm{ST}}(D, f), \eta) = p^{-1}|\hat{\xi}_{\mathrm{ST}}(D, f) - \xi|^2$ for $f \in \mathcal{F}_{\mathrm{ST}}$. The $i$th component of $\hat{\xi}_{\mathrm{ST}}(D, f)$ is

(4.23) $\quad \hat{\xi}_{\mathrm{ST},i}(D, f) = \mathrm{sgn}(z_i)(|z_i| - t)_+ = z_i - (|z_i| \wedge t)\,\mathrm{sgn}(z_i).$

Hence,

(4.24)
$$V(f) = p^{-1}\sum_{i=1}^p (z_i - \xi_i)^2 + p^{-1}\sum_{i=1}^p (|z_i| \wedge t)^2$$
$$- 2\sum_{i=1}^p (z_i - \xi_i)(|z_i| \wedge t)\,\mathrm{sgn}(z_i).$$

On the right-hand side of this equation, the $L_1$ convergence, uniformly over $t \geq 0$, of the second term is given by (4.9) and is immediate for the first term. It remains to verify this mode of convergence for

$$\sum_{i=1}^{p}(z_i - \xi_i)(|z_i| \wedge t)\,\mathrm{sgn}(z_i)$$

(4.25)
$$= \sum_{i=1}^{p}(z_i - \xi_i)z_i I(|z_i| \leq t) + \sum_{i=1}^{p}(z_i - \xi_i)t I(|z_i| > t)$$

$$= T_1(t) + T_2(t), \quad \text{say.}$$

For $i = 1, 2$, the analysis of $T_i(t)$ parallels that given for $S_i(t)$ after (4.13). Limit (4.2) now follows for $V(f) = L(\hat{\eta}_{\mathrm{ST}}(D, f), \eta)$.

(b) and (c) In analogy to $\hat{f} = \arg\min_{f \in \mathcal{F}} \hat{r}_{\mathcal{F}}(D, f)$, let $\tilde{f} = \arg\min_{f \in \mathcal{F}} r_{\mathcal{F}}(f, \xi, \sigma^2)$. Then $\min_{f \in \mathcal{F}} R(\hat{\eta}_{\mathcal{F}}(D, f), \eta, \sigma^2) = r_{\mathcal{F}}(\tilde{f}, \xi, \sigma^2)$. We first show that (4.2) implies

(4.26)
$$\lim_{p \to \infty} \sup_{\|\xi\| \leq c} \mathrm{E}|T - r_{\mathcal{F}}(\tilde{f}, \xi, \sigma^2)| = 0,$$

where $T$ can be $L(\hat{\eta}_{\mathcal{F}}(D, \hat{f}), \eta)$ or $L(\hat{\eta}_{\mathcal{F}}(D, \tilde{f}), \eta)$ or $\hat{r}_{\mathcal{F}}(D, \hat{f})$.

Indeed, (4.2) with $V(f) = \hat{r}_{\mathcal{F}}(D, f)$ entails

(4.27)
$$\lim_{p \to \infty} \sup_{\|\xi\| \leq c} \mathrm{E}|\hat{r}_{\mathcal{F}}(D, \hat{f}) - r_{\mathcal{F}}(\tilde{f}, \xi, \sigma^2)| = 0,$$

$$\lim_{p \to \infty} \sup_{\|\xi\| \leq c} \mathrm{E}|\hat{r}_{\mathcal{F}}(D, \hat{f}) - r_{\mathcal{F}}(\hat{f}, \xi, \sigma^2)| = 0.$$

Hence, (4.26) holds for $T = \hat{r}_{\mathcal{F}}(D, \hat{f})$ and

(4.28)
$$\lim_{p \to \infty} \sup_{\|\xi\| \leq c} \mathrm{E}|r_{\mathcal{F}}(\hat{f}, \xi, \sigma^2) - r_{\mathcal{F}}(\tilde{f}, \xi, \sigma^2)| = 0.$$

On the other hand, (4.2) with $V(f) = L(\hat{\eta}_{\mathcal{F}}(D, f), \eta)$ gives

(4.29)
$$\lim_{p \to \infty} \sup_{\|\xi\| \leq c} \mathrm{E}|L(\hat{\eta}_{\mathcal{F}}(D, \hat{f}), \eta) - r_{\mathcal{F}}(\hat{f}, \xi, \sigma^2)| = 0,$$

$$\lim_{p \to \infty} \sup_{\|\xi\| \leq c} \mathrm{E}|L(\hat{\eta}_{\mathcal{F}}(D, \tilde{f}), \eta) - r_{\mathcal{F}}(\tilde{f}, \xi, \sigma^2)| = 0.$$

These limits, together with (4.28), establish the remaining two cases of (4.26).

The limits (4.3) and (4.4) are immediate consequences of (4.26).

(d) By Theorem 2.1 of Beran and Dümbgen (1998), limit (4.2) with $\mathcal{F} = \mathcal{F}_{\mathrm{MS}}$ can be strengthened to

(4.30) $\displaystyle\sup_{\|\xi\| \leq c} \mathrm{E} \sup_{f \in \mathcal{F}_{\mathrm{MS}}} |V(f) - R(\hat{\eta}_{\mathcal{F}}(D, f), \eta, \sigma^2)| \leq C_1 p^{-1/2} + C_2 \mathrm{E}|\hat{\sigma}^2 - \sigma^2|,$

where the $C_i$ are finite constants. The first assertion of part (d) follows.

The arguments above for $\mathcal{F} = \mathcal{F}_{\mathrm{ST}}$ imply that

$$
\begin{aligned}
(4.31) \quad & \sup_{\|\xi\| \leq c} \mathrm{E} \sup_{f \in \mathcal{F}_{\mathrm{ST}}} |V(f) - R(\hat{\eta}_{\mathcal{F}}(D,f), \eta, \sigma^2)| \\
& \leq C_1 p^{-1/2} (\log(p))^{1/4} + C_2 \mathrm{E} |\hat{\sigma}^2 - \sigma^2|,
\end{aligned}
$$

where the $C_i$ are finite constants. The second assertion of part (d) follows.

(e) Part (e) similarly follows from (4.30) and (4.31).

Department of Statistics  
University of California  
Davis, California 95616  
USA  
e-mail: beran@wald.ucdavis.edu